\documentclass[12pt]{amsart}

\usepackage{amscd, amsmath, amssymb, amsthm, amsfonts, amsxtra, 
enumerate, latexsym, tabularx, setspace}
\usepackage[all]{xy}
\usepackage[dvips]{graphicx,color}
\usepackage{here}

\def\NN{{\mathbb{N}}}
\def\PP{{\mathbb{P}}}
\def\R{{\mathbb{R}}}
\def\ZZ{{\mathbb{Z}}}
\def\Proj{{\mathrm{Proj\; }}}
\def\Spec{{\mathrm{Spec\; }}}
\def\QQ{{\mathbb{Q}}}
\def\vol{{\mathrm{Vol}}}
\def\area{{\mathrm{Area}}}

\theoremstyle{plain}
\newtheorem{thm}{Theorem}[section]

\theoremstyle{definition}

\newtheorem{exa}[thm]{Example}

\theoremstyle{remark}

\title{Korff $F$-signatures of Hirzebruch surfaces}
\author{Daisuke Hirose}
\address{Department of General education, National Institute of Technology, Fukushima College, 
30 Aza-Nagao, Kamiarakawa, Iwaki-shi, Fukushima 970-8034, Japan}
\email{hirose@fukushima-nct.ac.jp}
\author{Tadakazu Sawada}
\address{Department of General Education, National Institute of Technology, Fukushima College, 
Fukushima 970-8034, Japan.}
\email{sawada@fukushima-nct.ac.jp}

\begin{document}
\maketitle
\markboth{Daisuke Hirose and Tadakazu Sawada}{Korff $F$-signature of Hirzebruch surfaces}

\begin{abstract}
M. Von~Korff introduced the $F$-signature of a normal projective variety, and computed 
the $F$-signature of the product of two projective lines $\PP^1\times \PP^1$. 
In this paper, we compute $F$-signatures of Hirzebruch surfaces. 
\end{abstract}

\section*{Introduction}
Let $R$ be a $d$-dimensional Noetherian local ring of prime characteristic $p$ with 
perfect residue field. Let $F:R\rightarrow R$ be the Frobenius map, that is, $F(x)=x^p$ 
for all elements $x$ of $R$. We denote by $F^e_{\ast}R$ the $R$-module whose abelian 
group structure is inherited by $R$, and whose $R$-module structure is given by the 
$e$-th iterated Frobenius map $F^e$. If $F_{\ast}R$ is finitely generated as an 
$R$-module, we say that $R$ is {\it $F$-finite}. Let the $R$-module $F^e_{\ast}R$ 
be decomposed as $F^e_{\ast}R\cong R^{\oplus a_e}\oplus M_e$, where $M_e$ does 
not have free summands. In \cite{HL}, C. Huneke and G. Leuschke introduced the notion 
of the $F$-signature $s(R)$ of $R$ as 
$$s(R)=\lim_{e\rightarrow \infty}\dfrac{a_e}{p^{ed}}.$$
The existence of the limit $\lim_{e\rightarrow \infty}({a_e}/{p^{ed}})$ is not trivial. In \cite{T}, 
K. Tucker proved that the limit exists (under the assumption that $R$ is $F$-finite). The 
notion of the $F$-signature is generalized by M. Hashimoto and Y. Nakajima in \cite{HN}, and 
A. Sannai in \cite{San}. As a global analogue, M. Von~Korff introduced the $F$-signature of a 
normal projective variety, which we call the Korff $F$-signature, in \cite{K}. He computed 
the Korff $F$-signature of the product of two projective lines $\PP^1\times \PP^1$. In this 
paper, we compute Korff $F$-signatures of Hirzebruch surfaces. 

In Section 1, we review generalities on Korff $F$-signatures. In Section 2, we compute Korff 
$F$-signatures of projective spaces and Hirzebruch surfaces. In Section 3 and 4, we gather 
together the figures and the source code of wxMaxima in the computation of $F$-signatures 
of Hirzebruch surfaces. 

\section{Preliminaries}
In this section, we review generalities on Korff $F$-signatures. See \cite{K} for details. 

Let $k$ be a field of positive characteristic $p$. Let $X$ be a normal projective variety over 
$k$, and let $D$ be a $\QQ$-divisor on $X$. We denote by $\mathrm{Sec}(X, D)$ the section 
ring $\bigoplus_{n\geq 0}\Gamma (X,\mathcal{O}_X(nD))$ of $D$. Suppose that 
$\mathrm{Sec}(X, D)$ is a finitely generated $k$-algebra of dimension at least two. If $c \in \NN$ 
is sufficiently divisible, then there exists an ample line bundle $L$ on $X$ such that 
$\mathrm{Sec}(X, cD)$ is isomorphic to a normal section ring 
$\mathrm{Sec}(X, L)=\bigoplus_{n\geq 0}\Gamma (X,L^{\otimes n})$ of $L$. We assume that 
$c$ is sufficiently divisible, and define the {\it Korff $F$-signature $s(X,D)$ of $X$ along $D$} 
to be $c\cdot s(\mathrm{Sec}(X, cD))$. (The $F$-signature of a graded ring is defined in the 
same way as the case of local rings.) 

Let $N\cong \ZZ^n$ be a lattice, and let $M$ be the dual lattice of $N$. We denote 
$N\otimes_{\ZZ} \R$ (resp. $M\otimes_{\ZZ} \R$) by $N_{\R}$ (resp. $M_{\R}$). Let $\Sigma$ 
be a fan in $N_{\R}$. We denote by $X_{\Sigma}$ the toric variety corresponding to $\Sigma$. 

Let $\Sigma$ be a complete fan in $N_{\R}$ with the primitive generators $v_1,\ldots ,v_n$. Let 
$D_1,\ldots, D_n$ be the torus-invariant prime Weil divisors corresponding to $v_1,\ldots ,v_n$, 
respectively. Let $D=\sum_{i=1}^na_iD_i$ be a torus-invariant Weil divisor on $X$. The 
{\it polytope associated to the divisor $D$} is defined to be the polytope 
$\{u\in M_{\R}|u\cdot v_i\geq -a_i {\rm\ for\ all}\ i\}$, and denoted by $P_D$. We may assume 
that $P_D$ contains the origin. Then we denote by $M_D$ the lattice $M\cap (P_D\cdot \R)$, 
where $P_D\cdot \R$ is the $\R$-vector subspace spanned by $P_D$ in $M_{\R}$. We denote 
by $N_D$ the dual lattice of $M_D$. We define $I_D$ to be the subset of $\{1,\ldots ,n\}$ such 
that the hyperplane $\{u|u\cdot v_i=-a_i\}$ determines a facet of $P_D$ for $i\in I_D$, and 
$P_{\Sigma, D}$ to be the polytope 
$\{(v, t)\in (M_D\times \ZZ)_{\R}| 0\leq (v,t)\cdot (v_i,a_i)<1 {\rm\ for\ all}\ i\in I_D\}$. 

Suppose that $P_D$ is not full-dimensional. Let $\pi: N\rightarrow N_D$ be the natural projection 
map. Let $c_i$ be the rational number such that $\pi (c_i v_i)$ is the primitive generator for its 
ray in $(N_D)_{\R}$. We define ${P'}_{\Sigma, D}$ to be the polytope 
$\{(v, t)\in (M_D\times \ZZ)_{\R}|0\leq (v,t)\cdot c_i(v_i,a_i)<1 {\rm\ for\ all}\ i\in I_D \}$. 

\begin{thm}[\cite{K}, Proposition 4.4.4, Corollary 4.4.7]
Let $X$ be a $d$-dimensional projective toric variety over $k$, and let $\Sigma$ be the fan 
corresponding to $X$. Let $v_1,\ldots ,v_n$ be the primitive generators of $\Sigma$, and let 
$D_1,\ldots, D_n$ be the torus-invariant prime Weil divisors corresponding to $v_1,\ldots ,v_n$, 
respectively. Let $D=\sum_{i=1}^na_iD_i$ be a (non-zero) effective torus-invariant Weil divisor 
on $X$. Then:\\
{\rm (1)} If $P_D$ is full-dimensional, then $s (X, D)=\vol(P_{\Sigma, D})$. \\
{\rm (2)} If $P_D$ is not full-dimensional, then $s (X, D)=\vol({P'}_{\Sigma, D})$. 
\end{thm}

\section{Examples}
M. Von~Korff computed Korff $F$-signatures of some toric varieties including the product of 
two projective lines $\PP^1\times \PP^1$ in \cite{K}. In this section, we compute Korff 
$F$-signatures of projective spaces and Hirzebruch surfaces. 
\begin{exa}[Projective spaces $\PP^n$]
Let $N=\ZZ^n$ with standard basis $e_1, \ldots ,e_n$. Let $e_0=-\sum_{i=1}^n e_i$, and let 
$\Sigma$ be the fan in $N_{\R}$ consisting of the cones generated by all proper subsets of 
$\{e_0,e_1,\ldots, e_n\}$. We see that $X_{\Sigma}$ is isomorphic to the $n$-dimensional 
projective space $\PP^n$. Let $D_i$ be the torus-invariant prime Weil divisor corresponding 
to the ray $\R_{+} e_i$ in $N_{\R}$ for $0\leq i \leq n$. 

Let $D=\sum_{i=0}^n a_i D_i$ be  a (non-zero) effective torus-invariant Weil divisor. We see 
that 
$$P_D=\left\{(x_1, \ldots , x_n)\in M_{\R} \left|\begin{array}{l}
\sum_{i=1}^n x_i \leq a_0, \\
x_1 \geq -a_1, \\
\cdots \\
x_n \geq -a_n 
\end{array}
\right.\right\}.$$
Since $a_i\not= 0$ for some $i$, $P_D$ has an interior point, and the equations 
$\sum_{i=1}^n x_i =a_0, x_1=-a_1,\ldots , x_n=-a_n$ define the facets of $P_D$, respectively. 
Hence we have $I_D=\{0, 1, \ldots ,n\}$. In what follows, we compute the volume of the 
$n$-dimensional parallelepipeds
$$P_{\Sigma, D}=
\left\{(x_1, \ldots , x_n, t)\in (M_D\times \ZZ)_{\R}\left|\begin{array}{l}
0\leq -\sum_{i=1}^n x_i+a_0 t <1, \\
0\leq x_1+a_1t <1, \\
\cdots \\
0\leq x_n+a_nt <1
\end{array}
\right.\right\}.$$

Let $a=\sum_{i=0}^{n} a_i$. The point of the intersection of the $n+1$ hyperplanes defined 
by $-\sum_{i=1}^n x_i+a_0t =1, x_1+a_1t =0,\ldots, x_n+a_nt =0$ is $(-a_1/a, \ldots ,-a_n/a, 1/a)$. 
This is given by the following elementary transformations of matrices:\\
$\left(\begin{array}{ccccc|c}
-1 & -1 & \cdots & -1 & a_0 & 1 \\
1 & 0 & \cdots & 0 & a_1 & 0\\
  &   & \cdots &    &      & \\
0 & 0 & \cdots & 1 & a_n & 0 
\end{array}
\right)
\rightarrow \left(\begin{array}{ccccc|c}
0 & 0 & \cdots & 0 & a & 1\\
1 & 0 & \cdots & 0 & a_1 & 0 \\
  &   &  \cdots &    &      & \\
0 & 0 & \cdots & 1 & a_n & 0 
\end{array}
\right)\\
\rightarrow \left(\begin{array}{ccccc|c}
0 & 0 & \cdots & 0 & a & 1 \\
1 & 0 &  \cdots & 0 & 0 & -a_1/a \\
  &   &  \cdots &    &      & \\
0 & 0 & \cdots & 1 & 0 & -a_n/a
\end{array}
\right)
\rightarrow \left(\begin{array}{ccccc|c}
1 & 0 & \cdots & 0 & 0 & -a_1/a \\
  &   & \cdots &    &      & \\
0 & 0 & \cdots & 1 & 0 & -a_n/a \\
0 & 0 & \cdots & 0 & 1 & 1/a
\end{array}
\right)$. 

By the same argument, we see that the point of the intersection of the $n+1$ hyperplanes 
defined by $-\sum_{i=1}^n x_i+a_0t =0$, $x_l+a_l t =0\ (l\not= k)$, $x_k+a_k t =1$ is 
$(-a_1/a,\ldots ,1-a_k/a,\ldots ,-a_n/a,1/a)$ for each $1\leq k\leq n$. 

Let $v_0={}^t(-a_1/a, \ldots ,-a_n/a, 1/a)$, and 
$$v_k={}^t(-a_1/a,\ldots ,1-a_k/a,\ldots ,-a_n/a,1/a)$$
for $1\leq k\leq n$. Since $P_{\Sigma, D}$ 
is spanned by $v_0, v_1,\ldots ,v_n$, we have 
\begin{eqnarray*}
s(\PP^n,D)&=&\vol (P_{\Sigma, D})\\
&=& |\mathrm{det}(v_0, v_1, \ldots ,v_n)|\\
&=& \left|\mathrm{det}\left(\begin{array}{ccccc}
-a_1/a & 1-a_1/a & -a_1/a & \cdots &  -a_1/a \\
-a_2/a & -a_2/a & 1-a_2/a & \cdots &  -a_2/a \\
 &   &  & \cdots & \\
-a_n/a & -a_n/a & -a_n/a & \cdots & 1-a_n/a \\
1/a & 1/a & 1/a & \cdots & 1/a 
\end{array}\right)
\right|\\
&=&1/a^{n+1}
\left|\mathrm{det}\left(\begin{array}{ccccccc}
-a_1 & a-a_1 & -a_1 & & \cdots & &  -a_1 \\
-a_2 & -a_2 & a-a_2 & & \cdots & &  -a_2 \\
  &   &    &     \cdots &    &      & \\
-a_n & -a_n & -a_n &  & \cdots & & a-a_n \\
1 &  1 & 1 & \cdots & &  1 & 1 
\end{array}\right)
\right|\\
&=&1/a^{n+1}
\left|\mathrm{det}\left(\begin{array}{cccccc}
0 & a &  0 & & \cdots  & 0 \\
0 & 0 &  a & & \cdots  & 0 \\
  &   &    &    & \cdots         & \\
0 & 0 &  0 &  & \cdots  & a \\
 1 &  1 & 1 &  & \cdots  & 1 
\end{array}\right)
\right|\\
&=&1/a=1/(a_0+a_1+\cdots +a_n).
\end{eqnarray*}
\end{exa}

The following computation is the main result. 

\begin{exa}[Hirzebruch surfaces]
	All figures of this example are in Section \ref{fig}.
	Let $r$ be a positive integer, and let $N=\ZZ^2$ with standard basis $e_1$ and $e_2$.
	Let $v_1=e_1, v_2=e_2, v_3=-e_2$, and $v_4=-e_1+re_2$.
	We consider a fan $\Sigma_r$ in $N_{\R}=\R^2$.
	The fan $\Sigma_r$ consists the four two-dimensional cones $\langle v_1, v_2\rangle, \langle v_1, v_3\rangle, \langle v_2, v_4\rangle$, and $\langle v_3, v_4\rangle$, and all their faces.
	The Hirzebruch surface $X=X_{\Sigma_r}$ is defined by the fan $\Sigma_r$
(see, e.g., \cite{CLS}, Example 3.1.16.).
	Let $D_i$ be the torus-invariant prime divisors of $X$ corresponding to $v_i$.
	We have linear equivalents $D_1 \sim D_4$ and $D_3 \sim D_2+rD_4$ in the divisor class group of $X$.
	For an effective $\QQ$-divisor $D=\sum_{i=1}^4 a_i D_i$, we have
		\begin{eqnarray*}
			D\sim & a_1D_1+(a_2+a_3)D_2+(ra_3+a_4)D_4\\
		 	 \sim &(a_1+ra_3+a_4)D_1+(a_2+a_3)D_2.
		\end{eqnarray*}
	Hence $\mathcal{O}(D)\cong\mathcal{O}((a_1+ra_3+a_4)D_1+(a_2+a_3)D_2)$.
	Without loss of generality, we only compute the $F$-signature $s(X,D)$ of $X$ with respect to a (non-zero) effective $\QQ$-divisor $D=a_1D_1+a_2D_2$.

The polytope $P_D$ is to be
	\begin{eqnarray*}
		P_D=\left\{(x,y)\in M_{\R}\left|\begin{array}{l}
				x\geq -a_1,\ y\geq -a_2\\
				y\leq 0,\ -x+ry\geq 0\end{array}
					\right\}\right..
	\end{eqnarray*}
	
By the figure \ref{fig001}, we see that
$P_D$ is a square if and only if $ra_2<a_1$.
In that case, $I_D=\{1,2,3,4\}$.
On the other hand, if $ra_2\geq a_1$, then $P_D$ is a triangle.
That implies $I_D=\{1,3,4\}$.

\noindent(1) Suppose that $P_D$ is a square, i.e.,
$ra_2<a_1$.
Then $I_D=\{1,2,3,4\}$
Therefore
	\begin{eqnarray*}
		P_{\Sigma,D}= &\left\{(x,y,t)\in (M_D\times \ZZ)_{\R}\left|
			\begin{array}{l}
				0\leq x+a_1 t < 1,\\
				0\leq y+a_2 t < 1,\\
				0\leq -y < 1,\\
				0\leq -x+ry < 1
			\end{array} \right.\right\}.
	\end{eqnarray*}
We define polygons $Q$, $R(x)$, and $S(x)$ to be 
	\begin{eqnarray*}
		Q&=&\left\{(y,t)\in \R^2\left|
			\begin{array}{l}
				0\leq y+a_2 t <1,\\
				0\leq -y < 1
			\end{array} \right.\right\}\\
		 &=&\left\{(y,t)\in \R^2\left|
			\begin{array}{l}
				-a_2 t\leq y < 1-a_2 t,\\
				-1 < y\leq 0
			\end{array} \right.\right\},\\
	R(x)&=&\left\{(y,t)\in \R^2\left|
			\begin{array}{l}
				0\leq x+a_1 t <1,\\
				0\leq -x+ry < 1
			\end{array} \right.\right\}\\
	    &=&\left\{(y,t)\in \R^2\left|
			\begin{array}{l}
				-x/a_1\leq t < -x/a_1+1/a_1,\\
				x/r \leq y < x/r+1/r
			\end{array} \right.\right\},\\
	\end{eqnarray*}
and $S(x)=Q\cap R(x)$.(See the figure \ref{fig002}.)

Then $Q$ is a parallelogram, $R(x)$ is a rectangle, and
	\begin{eqnarray*}
		P_{\Sigma,D}=\{(x,y,t)\in (M_D\times \ZZ)_{\R}|x_{\min}\le x\le x_{\max},(y,t)\in S(x)\},
	\end{eqnarray*}
where $x_{\min}=\min\{x\in\R|S(x)\neq\emptyset \}$
and $x_{\max}=\max\{x\in\R|S(x)\neq\emptyset\}$.
Therefore
	\begin{eqnarray*}
		s(X,D)=\mathrm{Vol}(P_{\Sigma,D})=\int_{x_{\min}}^{x_{\max}} \area(S(x))\,dx.
	\end{eqnarray*}
	
We define some points and lines in the $ty$-plane appearing in our argument.	
Let $A,B,C$, and $D$ be the vertices of the rectangle $R(x)$.
In particular,
		\begin{eqnarray*}
		A=\left(-\dfrac{1}{a_1}x,\dfrac{1}{r}x\right),
		B=\left(-\dfrac{1}{a_1}x+\dfrac{1}{a_1},\dfrac{1}{r}x\right),\\
		C=\left(-\dfrac{1}{a_1}x+\dfrac{1}{a_1},\dfrac{1}{r}x+\dfrac{1}{r}\right),
		\mbox{and }
		D=\left(-\dfrac{1}{a_1}x,\dfrac{1}{r}x+\dfrac{1}{r}\right).
	\end{eqnarray*}
The point $A$ moves along the line $y=-a_1t/r$, denoted by $l_1$.
The points $B$ and $D$ move along the line $y=-a_1t/r+1/r$, denoted by $l_2$.
The point $C$ moves along the line $y=-a_1t/r+2/r$, denoted by $l_3$.

We denote the vetices of the parallelogram $Q$
except the origin $O$ by $E,F$, and $G$.
That is,
		\begin{eqnarray*}
			E=\left(\dfrac{1}{a_2},-1\right),
			F=\left(\dfrac{2}{a_2},-1\right),
			\mbox{and }
			G=\left(\dfrac{1}{a_2},0\right).
		\end{eqnarray*}
Let $H$ and $I$ be the intersection points of the $t$-axis 
with $l_2$ and $l_3$, respectively.
We denote the intersection points of the line $y=-1$ 
with $l_1$, $l_2$, and $l_3$ by $J$, $K$, and $L$, respectively.
That is,
		\begin{eqnarray*}
			H=\left(\dfrac{1}{a_1},0\right),&
			I=\left(\dfrac{2}{a_1},0\right),&\\
			J=\left(\dfrac{r}{a_1},-1\right),&
			K=\left(\dfrac{r+1}{a_1},-1\right),&
			\mbox{and }
			L=\left(\dfrac{r+2}{a_1},-1\right).
		\end{eqnarray*}

For the point $A$, we denote the $t$-coordinate of $A$ and $y$-coordinate of $A$ by $A_t$ and $A_y$, respectively.
We use the same notation for all points in the $ty$-plane.
For example, $B_t=-x/a_1+1/a_1$ and $J_y=-1$.

The shape of the polygon $P_{\Sigma, D}$ depends on how the lines $l_1, l_2$ and $l_3$ across the parallelogram $Q$.
We divide into five subcases from (1-1) through (1-5):

\noindent(1-1) 
We assume that $L_t\leq E_t$.
That is, $(r+2)/a_1\le 1/a_2$.
This is equivalent to $(r+2)a_2\le a_1$.
Then $2a_2\le a_1$.
This is equivalent to $2/a_1\le 1/a_2$.
Hence $I_t\le G_t$.
(See the figure \ref{fig101}.)

We denote the intersection points of the line $OE$ with $l_2$ and $l_3$ by $T$ and $U$, respectively.
Then
	\begin{eqnarray*}
		T=\left(\frac{1}{a_1-ra_2},-\frac{a_2}{a_1-ra_2}\right),\ \mbox{and }
		U=\left(\frac{2}{a_1-ra_2},-\frac{2a_2}{a_1-ra_2}\right).
	\end{eqnarray*}

If $C=U$, then $(-x/a_1+1/a_1,x/r+1/r)=(2/(a_1-ra_2),-2a_2/(a_1-ra_2))$.
Hence $x=-(a_1+ra_2)/(a_1-ra_2)$.
If $D=T$, then $(-x/a_1,x/r+1/r)=(1/(a_1-ra_2),-a_2/(a_1-ra_2))$.
Hence $x=-a_1/(a_1-ra_2)$.
If $B=T$, then $(-x/a_1+1/a_1,x/r)=(1/(a_1-ra_2),-a_2/(a_1-ra_2))$.
Hence $x=-ra_2/(a_1-ra_2)$.
If $C=I$ and $D=H$, then $(-x/a_1+1/a_1,x/r+1/r)=(2/a_1,0)$.
Hence $x=-1$.
If $A=O$ and $B=H$, then $(-x/a_1,x/r)=(0,0)$.
Hence $x=0$.

Since $0<ra_2<a_1$, we have $ra_2<a_1<a_1+ra_2$ and $a_1-ra_2>0$.
Hence $-(a_1+ra_2)/(a_1-ra_2)<-a_1/(a_1-ra_2)<-ra_2/(a_1-ra_2)$.
Since $a_1-ra_2<a_1$, we have $-a_1/(a_1-ra_2)<-1$.
If $2ra_2\leq a_1$,
then $-1\leq -ra_2/(a_1-ra_2)$.
Therefore
	\[
		-\frac{a_2+ra_2}{a_1-ra_2}<-\frac{a_1}{a_1-ra_2}<-1<-\frac{ra_2}{a_1-ra_2}<0.
	\]
On the other hand, if $2ra_2>a_1$,
then $-ra_2/(a_1-ra_2)<-1$.
Therefore
	\[
		-\frac{a_2+ra_2}{a_1-ra_2}<-\frac{a_1}{a_1-ra_2}<-\frac{ra_2}{a_1-ra_2}<-1<0.
	\]

For various values of $x$,
we consider shapes of $S(x)$ and areas of them.

If $-(a_2+ra_2)/(a_1-ra_2)<x<-a_1/(a_1-ra_2)$,
then $S(x)$ is a triangle as in the figure \ref{fig102}.
Hence
	\begin{eqnarray*}
		\area(S(x))&=&\frac{1}{2}\left\{\left(\frac{1}{r}x+\frac{1}{r}\right)-\left(\frac{a_2}{a_1}x-\frac{a_2}{a_1}\right)\right\}\cdot\left\{\left(-\frac{1}{a_1}x+\frac{1}{a_1}\right)-\left(-\frac{1}{ra_2}x-\frac{1}{ra_2}\right) \right\}\\
		&=&\frac{1}{2a_1^2a_2r^2}\left\{(a_1-ra_2)x+(a_1+ra_2)\right\}^2.
	\end{eqnarray*}

\noindent(1-1-1) Suppose that $2ra_2\le a_1$.

If $-a_1/(a_1-ra_2)<x<-1$,
then $S(x)$ is a trapezoid as in the figure \ref{fig111}.
Hence
	\begin{eqnarray*}
	\area(S(x))&=&\frac{1}{2}\left\{\left(\frac{1}{r}x+\frac{1}{r}-\frac{a_2}{a_1}x\right)+\left(\frac{1}{r}x+\frac{1}{r}-\frac{a_2}{a_1}x+\frac{a_2}{a_1}\right)\right\}\cdot\frac{1}{a_1}\\
	&=&\frac{1}{2a_1^2r}\left\{2(a_1-ra_2)x+2a_1+ra_2\right\}.
	\end{eqnarray*}
We denote this area by $V_1$.

If $-1<x<-ra_2/(a_1-ra_2)$,
then $S(x)$ is a tropezoid as in the figure \ref{fig112}.
Hence
	\begin{eqnarray*}
		\area(S(x))&=&\frac{1}{2}\left(-\frac{a_2}{a_1}x-\frac{a_2}{a_1}x+\frac{a_2}{a_1}\right)\cdot\frac{1}{a_1}\\
		&=&\frac{a_2}{2a_1^2}(1-2x).
	\end{eqnarray*}
We denote this area by $V_2$.

If $-ra_2/(a_1-ra_2)<x<0$,
then $S(x)$ is a pentagon as in the figure \ref{fig113}.
Hence
	\begin{eqnarray*}
		\area(S(x))&=&-\frac{1}{r}x\cdot\frac{1}{a_1}-\frac{1}{2}\left(\frac{a_2}{a_1}x-\frac{1}{r}x\right)\cdot\left\{-\frac{1}{ra_2}x-\left(-\frac{1}{a_1}x\right)\right\}\\
		&=&-\frac{1}{ra_1}x-\frac{(a_1-ra_2)^2}{2r^2a_1^2a_2}x^2.
	\end{eqnarray*}
We denote this area by $V_3$.

Therefore
	\begin{eqnarray*}
		s(X,D)&=&\vol(P_{\Sigma,D})\\
		&=&\int_{-\frac{a_2+ra_2}{a_1-ra_2}}^{-\frac{a_1}{a_1-ra_2}}\frac{1}{2a_1^2a_2r^2}\left\{(a_1-ra_2)x+(a_1+ra_2)\right\}^2dx\\
		&&+\int_{-\frac{a_1}{a_1-ra_2}}^{-1}V_1\, dx
		+\int_{-1}^{-\frac{ra_2}{a_1-ra_2}}V_2\, dx
		+\int_{-\frac{ra_2}{a_1-ra_2}}^0V_3\, dx\\
		&=&\frac{a_2}{a_1(a_1-a_2r)}.
	\end{eqnarray*}
In the above calculation, we use wxMaxima
(see Section \ref{maxima} (\%o1)).

\noindent(1-1-2)
Suppose that $2ra_2>a_1$.
If $-a_1/(a_1-ra_2)<x<-ra_2/(a_1-ra_2)$,
then $S(x)$ is a trapezoid as in the figure \ref{fig121},
which is the same shape as that in the case where 
 $-a_1/(a_1-ra_2)<x<-1$ on (1-1-1).
Then $\area(S(x))=V_1$.

If $-ra_2/(a_1-ra_2)<x<-1$,
then $S(x)$ is a pentagon as in the figure \ref{fig122}.
Hence
	\begin{eqnarray*}
		\area(S(x))&=&\left(\frac{1}{r}x+\frac{1}{r}-\frac{1}{r}\right)\cdot\frac{1}{a_1}
			-\frac{1}{2}\left(\frac{a_2}{a_1}x-\frac{1}{r}x\right)\cdot\left\{-\frac{1}{ra_2}x-\left(-\frac{1}{a_1}x\right)\right\}\\
			&=&\frac{1}{ra_1}-\frac{(a_1-ra_2)^2}{2a_1^2a_2r^2}x^2.
	\end{eqnarray*}
We denote this area by $V_4$.

If $-1<x<0$,
then $S(x)$ is a pentagon as in the figure \ref{fig123},
which is the same shape as that in the case where 
 $-ra_2/(a_1-ra_2)<x<0$ on (1-1-1).
Then $\area(S(x))=V_3$.

Therefore
	\begin{eqnarray*}
		s(X,D)&=&\vol(P_{\Sigma,D})\\
		&=&\int_{-\frac{a_2+ra_2}{a_1-ra_2}}^{-\frac{a_1}{a_1-ra_2}}\frac{1}{2a_1^2a_2r^2}\left\{(a_1-ra_2)x+(a_1+ra_2)\right\}^2dx\\
		&&+\int_{-\frac{a_1}{a_1-ra_2}}^{-\frac{ra_2}{a_1-ra_2}}V_1\,dx
		+\int_{-\frac{ra_2}{a_1-ra_2}}^{-1}V_4\,dx
		+\int_{-1}^0V_3\,dx\\
		&=&\frac{a_2}{a_1(a_1-a_2r)},
	\end{eqnarray*}
(see Section \ref{maxima} (\%o2)).

\noindent(1-2) 
We assume that $K_t< E_t< L_t$.
That is, $(r+1)/a_1< 1/a_2< (r+2)/a_1$.
This is equivalent to $(r+1)a_2< a_1<(r+2)a_2$.
Then $2a_2< a_1$.
We have $I_t\le G_t$.
Since $(r+1)a_2<a_1$ and $a_2\leq ra_2<a_1$,
we have $(r+2)a_2<2a_1$.
Therefore $L_t<F_t$.
The point $T$ denotes the same in the case (1-1).
(See the figure \ref{fig201}.)

If $C=L$ and $D=K$,
then $(-x/a_1+1/a_1,x/r+1/r)=((r+2)/a_1,-1)$.
Hence $x=-(r+1)$.
If $B_t=E_t$,
then $-(1/a_1)x+1/a_1=1/a_2$.
Hence $x=-(a_1-a_2)/a_2$.
By the same argument in the case (1-1), 
we have the following four values of $x$.
If $D=T$, then $x=-a_1/(a_1-ra_2)$.
If $B=T$, then $x=-ra_2/(a_1-ra_2)$.
If $C=I$ and $D=H$, then $x=-1$.
If $A=O$ and $B=H$, then $x=0$.

Since $(r+1)a_2<a_1$,
we have $r(r+1)a_2+a_1<(r+1)a_1$.
Hence $a_1<(r+1)(a_1-ra_2)$.
Therefore $-(r+1)<-a_1/(a_1-ra_2)$.
If $a_1/(a_1-ra_2)>(a_1-a_2)/a_2$,
then $a_1^2-(r+2)a_2a_1+ra_2^2>0$.
We have 
	\begin{eqnarray*}
		\frac{(r+2)-\sqrt{r^2+4}}{2}a_2<a_1<\frac{(r+2)+\sqrt{r^2+4}}{2}a_2.
	\end{eqnarray*}
By the assumption that $(r+1)a_2<a_1<(r+2)a_2$,
we have $a_1/(a_1-ra_2)>(a_1-a_2)/a_2$ in case (1-2).
Therefore $-a_1/(a_1-ra_2)<-(a_1-a_2)/a_2$.
Since $2a_2<a_1$, we have $a_2<a_1-a_2$.
Hence $-(a_1-a_2)/a_2<-1$.
The assumption that $(r+1)a_2<a_1$ gives the inequalities $ra_2<a_1-a_2$ and $a_2<a_1-ra_2$.
Those imply $ra_2/(a_1-ra_2)<(a_1-a_2)/a_2$.
Therefore $-(a_1-a_2)/a_2<-ra_2/(a_1-ra_2)$.
Finally, we compare the inequality between $-1$ and $-ra_2/(a_1-ra_2)$.
Suppose that $-1<-ra_2/(a_1-ra_2)$.
Then $a_1>2ra_2$.
With the assumption that $a_1<(r+2)a_2$, we have $2ra_2<(r+2)a_2$.
Hence $r<2$.
That is, $r=1$.
Therefore if $r=1$,
then
	\[
	-(r+1)<-\frac{a_1}{a_1-ra_2}<-\frac{a_1-a_2}{a_2}<-1<-\frac{ra_2}{a_1-ra_2}<0.
	\]
By the same argument, if $-ra_2/(a_1-ra_2)<-1$, then $1<r$.
Therefore if $r\geq 2$,
then
	\[
	-(r+1)<-\frac{a_1}{a_1-ra_2}<-\frac{a_1-a_2}{a_2}<-\frac{ra_2}{a_1-ra_2}<-1<0.
	\]
We remark that $-ra_2/(a_1-ra_2)\neq -1$, because $r$ is integer.

For various values of $x$,
we consider shapes of $S(x)$ and areas of them.

If $-(r+1)<x<-a_1/(a_1-ra_2)$,
then $S(x)$ is a trapezoid as in the figure \ref{fig202}.
Hence
	\begin{eqnarray*}
		\area(S(x))&=&\frac{1}{2}\left\{\left(-\frac{1}{a_1}x+\frac{1}{a_1}-\left(-\frac{1}{ra_2}x-\frac{1}{ra_2}\right)\right)+\left(-\frac{1}{a_1}x+\frac{1}{a_1}-\frac{1}{a_2}\right)\right\}\\
		&&\times\left(\frac{1}{r}x+\frac{1}{r}-(-1)\right) \\
		&=&\frac{1}{2a_1a_2r^2}\left\{(a_1-2a_2r)x+2a_2r+a_1-a_1r\right\}(x+r+1).
	\end{eqnarray*}

If $-a_1/(a_1-ra_2)<x<-(a_1-a_2)/a_2$,
then $S(x)$ is a pentagon as in the figure \ref{fig203}.
Hence
	\begin{eqnarray*}
		\area(S(x))&=&\left\{\frac{1}{r}x+\frac{1}{r}-(-1)\right\}\cdot\frac{1}{a_1}
		-\frac{1}{2}\left\{\frac{1}{a_2}-\left(-\frac{1}{a_1}x\right)\right\}\cdot\left\{\frac{a_2}{a_1}x-(-1)\right\}\\
		&=&\frac{1}{a_1r}(x+r+1)-\frac{1}{2a_1^2a_2}(a_2x+a_1)^2.
	\end{eqnarray*}
We denote this area by $V_5$.

\noindent(1-2-1) Suppose that $r=1$.
If $-(a_1-a_2)/a_2<x<-1$,
then $S(x)$ is a trapezoid as in the figure \ref{fig211},
which is the same shape as that in the case where 
$-a_1/(a_1-ra_2)<x<-1$ on (1-1-1).
Hence $\area(S(x))=V_1$.

If $-1<x<-ra_2/(a_1-ra_2)$,
then $S(x)$ is a trapezoid as in the figure \ref{fig212},
which is the same shape as that in the case where 
$-1<x<-ra_2/(a_1-ra_2)$ on (1-1-1).
Hence $\area(S(x))=V_2$.

If $-ra_2/(a_1-ra_2)<x<0$,
then $S(x)$ is a pentagon as in the figure \ref{fig213},
which is the same shape as that in the case where 
$-ra_2/(a_1-ra_2)<x<0$ on (1-1-1).
Hence $\area(S(x))=V_3$.

Therefore
	\begin{eqnarray*}
		s(X,D)&=&\vol(P_{\Sigma,D})\\
		&=&\int_{-(r+1)}^{-\frac{a_1}{a_1-ra_2}}\frac{1}{2a_1a_2r^2}\left\{(a_1-2a_2r)x+2a_2r+a_1-a_1r\right\}(x+r+1)dx\\
		&&+\int_{-\frac{a_1}{a_1-ra_2}}^{-\frac{a_1-a_2}{a_2}}V_5\,dx
		+\int_{-\frac{a_1-a_2}{a_2}}^{-1}V_1\,dx
		+\int_{-1}^{-\frac{ra_2}{a_1-ra_2}}V_2\,dx
		+\int_{-\frac{ra_2}{a_1-ra_2}}^{0}V_3\,dx\\
		&=&\frac{1}{6a_1a_2^2(a_2r-a_1)}(a_2^3r^3+6a_2^3r^2-3a_1a_2^2r^2+12a_2^3r-12a_1a_2^2r\\
		&&+3a_1^2a_2r+2a_2^3-12a_1a_2^2+6a_1^2a_2-a_1^3),
	\end{eqnarray*}
(see Section \ref{maxima} (\%o3)).

\noindent(1-2-2) Suppose that $r\geq 2$.
If $-(a_1-a_2)/a_2<x<-ra_2/(a_1-ra_2)$,
then $S(x)$ is a trapezoid as in the figure \ref{fig221},
which is the same shape as that in the case where 
$-(a_1-a_2)/a_2<x<-1$ on (1-2-1).
Hence $\area(S(x))=V_1$.

If $-ra_2/(a_1-ra_2)<x<-1$,
then $S(x)$ is a pentagon as in the figure \ref{fig222},
which is the same shape as that in the case where 
$-ra_2/(a_1-ra_2)<x<-1$ on (1-1-2).
Hence $\area(S(x))=V_4$

If $-1<x<0$,
then $S(x)$ is a pentagonas as in the figure \ref{fig223},
which is the same shape as that in the case where 
$-ra_2/(a_1-ra_2)<x<0$ on (1-2-1).
Hence $\area(S(x))=V_3$.

Therefore
	\begin{eqnarray*}
		s(X,D)&=&\vol(P_{\Sigma,D})\\
		&=&\int_{-(r+1)}^{-\frac{a_1}{a_1-ra_2}}\frac{1}{2a_1a_2r^2}\left\{(a_1-2a_2r)x+2a_2r+a_1-a_1r\right\}(x+r+1)dx\\
		&&+\int_{-\frac{a_1}{a_1-ra_2}}^{-\frac{a_1-a_2}{a_2}}V_5\,dx
		+\int_{-\frac{a_1-a_2}{a_2}}^{-\frac{ra_2}{a_1-ra_2}}V_1\,dx
		+\int_{-\frac{ra_2}{a_1-ra_2}}^{-1}V_4\,dx
		+\int_{-1}^{0}V_3\,dx\\
		&=&\frac{1}{6a_1a_2^2(a_2r-a_1)}(
		a_2^3r^3+6a_2^3r^2-3a_1a_2^2r^2+12a_2^3r-12a_1a_2^2r\\
		&&+3a_1^2a_2r+2a_2^3-12a_1a_2^2+6a_1^2a_2-a_1^3),
	\end{eqnarray*}
(see Section \ref{maxima} (\%o4)).

\noindent(1-3) 
We assume that $E_t\leq K_t$.
That is, $1/a_2\leq (r+1)/a_1$.
This is equivalent to $a_1\leq(r+1)a_2$.
We also assume that $r\geq 2$.
Then $a_1>ra_2\geq 2a_2$.
Hence $2/a_1<1/a_2$ and $(r+2)/a_1<2/a_2$.
Hence $I_t<G_t$ and $L_t<F_t$.
(See the figure \ref{fig301}.)

By the same argument in the case (1-2), $C=L$ and $D=K$ imply $x=-(r+1)$.
If $A_t=E_t$,
then $-x/a_1=1/a_2$.
Hence $x=-a_1/a_2$.
If $B=K$, then $(-x/a_1+1/a_1,x/r)=((r+2)/a_1,-1)$.
Hence $x=-r$.
By the same argument in the case (1-1), we have the following two values of $x$.
If $C=I$ and $D=H$, then $x=-1$.
If $A=O$ and $B=H$, then $x=0$.

Since $ra_2<a_1\leq (r+1)a_2$, we have $-(r+1)\leq -a_1/a_2<-r$.
Therefore
	\[
	-(r+1)<-\frac{a_1}{a_2}<-r<-1<0.
	\]

For various values of $x$,
we consider shapes of $S(x)$ and areas of them.

If $-(r+1)<x<-a_1/a_2$,
then $S(x)$ is a rectangle as in the figure \ref{fig302}.
Hence
	\begin{eqnarray*}
		\area(S(x))&=&\left\{\frac{1}{r}x+\frac{1}{r}-(-1)\right\}\cdot\frac{1}{a_1}\\
		&=&\frac{1}{a_1r}(x+r+1).
	\end{eqnarray*}
We denote this area by $V_6$.

If $-a_1/a_2<x<-r$,
then $S(x)$ is a pentagon as in the figure \ref{fig303},
which is the same shape as that in the case where $-a_1/(a_1-ra_2)<x<-(a_1-a_2)/a_2$ on (1-2).
Hence $\area(S(x))=V_5$.

If $-r<x<-1$,
then $S(x)$ is a pentagon as in the figure \ref{fig304},
which is the same shape as that in the case where 
$-ra_2/(a_1-ra_2)<x<-1$ on (1-1-2).
Hence $\area(S(x))=V_4$.

If $-1<x<0$,
then $S(x)$ is a pentagon as in the figure \ref{fig305},
which is the same shape as that in the case where 
$-ra_2/(a_1-ra_2)<x<0$ on (1-2-1).
Hence $\area(S(x))=V_3$.

Therefore
	\begin{eqnarray*}
		s(X,D)&=&\vol(P_{\Sigma,D})\\
		&=&\int_{-(r+1)}^{-\frac{a_1}{a_2}}V_6\,dx
		+\int_{-\frac{a_1}{a_2}}^{-r}V_5\,dx
		+\int_{-r}^{-1}V_4\,dx
		+\int_{-1}^{0}V_3\,dx\\
		&=&-\frac{a_2^2r^2-2a_1a_2r-6a_2^2+a_1^2}{6a_1a_2^2},
	\end{eqnarray*}
(see Section \ref{maxima} (\%o5)).

\noindent(1-4) 
We assume that $E_t\leq K_t$.
By the same argument in the case (1-3), this is equivalent to $a_1\leq (r+1)a_2$.
We assume that $r=1$.
Then $a_1\leq 2a_2$, thus $1/a_2\leq 2/a_1$.
Hence $G_t\leq I_t$.
We also assume that $3a_2<2a_1$.
Then $3/a_1<2/a_2$.
Hence $L_t<F_t$.
(See the figure \ref{fig401}.)

We denote the intersection points of the line $FG$ with $l_3$ by $V$.
Then
	\begin{eqnarray*}
		V=\left(\frac{2-r}{a_1-ra_2},-\frac{2a_2-a_1}{a_1-ra_2}\right).
	\end{eqnarray*}
Since $r=1$, we have $V=(1/(a_1-a_2),(a_1-2a_2)/(a_1-a_2))$.

If $C=V$,
then we have $(-x/a_1+1/a_1,x/r+1/r)=((2-r)/(a_1-ra_2),-(a_1-2a_2)/(a_1-ra_2))$.
Hence $x=-(a_1+r(a_2-a_1))/(a_1-ra_2)$.
In particular, $x=-a_2/(a_1-a_2)$.
By the same argument in the case (1-3), we have the following two values of $x$.
If $A_t=E_t$, then $x=-a_1/a_2$.
If $B=K$, then $x=-r$, in particular $x=-1$.
By the same argument in the case (1-2), we have the following two values of $x$.
If $C=L$ and $D=K$,
then $x=-(r+1)$, in particular, $x=-2$.
If $B_t=E_t$, then $x=-(a_1-a_2)/a_2$.
By the same argument in the case (1-1), we have the following two values of $x$.
If $C=I$ and $D=H$,
then $x=-1$.
If $A=O$ and $B=H$,
then $x=0$.

Since $a_1\leq 2a_2$,
we have $-(r+1)=-2\leq -a_1/a_2$.
Since $3a_2<2a_1$, we have $a_2<2(a_1-a_2)$.
Hence $-2\leq -a_2/(a_1-a_2)$.
If $-a_2/(a_1-a_2)>-a_1/a_2$,
then $a_1^2-a_2a_1-a_2^2>0$.
Therefore $2a_1<(1-\sqrt{5})a_2$ or $(1+\sqrt{5})a_2< 2a_1$.
Since $0<a_2<a_1$, we have $2a_1>(1+\sqrt{5})a_2$.
On the other hand,
if $-a_2/(a_1-a_2)\leq -a_1/a_2$,
then $2a_1\leq (1+\sqrt{5})a_2$.
Since $a_2<a_1$, we have $-a_1/a_2<-1$.
Since $a_1\leq 2a_2$, we have $a_1-a_2\leq a_2$.
Hence $-a_2/(a_1-a_2)\leq -1 \leq -(a_1-a_2)/a_2$.

For various values of $x$,
we consider shapes of $S(x)$ and areas of them.

\noindent(1-4-1) Suppose that $2a_1>(1+\sqrt{5})a_2$.
Then we have
	\[
	-(r+1)=-2<-\frac{a_1}{a_2}<-\frac{a_2}{a_1-a_2}<-1<-\frac{a_1-a_2}{a_2}<0.
	\]
If $-(r+1)=-2<x<-a_1/a_2$,
then $S(x)$ is a rectangle as in the figure \ref{fig411},
which is the same shape as that in the case where $-(r+1)<x<-a_1/a_2$ on (1-3).
Hence $\area(S(x))=V_6$.

If $-a_1/a_2<x<-a_2/(a_1-a_2)$,
then $S(x)$ is a pentagon as in the figure \ref{fig412},
which is the same shape as that in the case where 
$-a_1/(a_1-ra_2)<x<-(a_1-a_2)/a_2$ on (1-2-1).
Hence $\area(S(x))=V_5$.

If $-a_2/(a_1-a_2)<x<-1$,
then $S(x)$ is a hexagon as in the figure \ref{fig413}.
Hence
	\begin{eqnarray*}
		\area(S(x))&=&\frac{1}{a_1}\left(\frac{1}{r}x+\frac{1}{r}+1\right)
				-\frac{1}{2}\left(\frac{a_2}{a_1}x+1\right)\cdot\left(\frac{1}{a_2}+\frac{1}{a_1}\right)\\
				&&-\frac{1}{2}\left\{-\frac{1}{a_1}x+\frac{1}{a_1}-\left(\frac{1}{a_2}-\frac{1}{ra_2}x-\frac{1}{ra_2}\right)\right\}\cdot\left\{\frac{1}{r}x+\frac{1}{r}-\left(\frac{a_2}{a_1}x-\frac{a_2}{a_1}+1\right)\right\}\\
				&=&\frac{1}{ra_1}(x+r+1)-\frac{1}{2a_1^2a_2}(a_2x+a_1)^2\\
				&&-\frac{1}{2a_1^2a_2r^2}\left\{(a_1-ra_2)x+a_1+ra_2-ra_1\right\}^2.
	\end{eqnarray*}
We denote this area by $V_7$.

If $-1<x<-(a_1-a_2)/a_2$,
then $S(x)$ is a hexagon as in the figure \ref{fig414}.
Hence
	\begin{eqnarray*}
		\area(S(x))&=&\frac{1}{a_1}\cdot\left(-\frac{1}{r}x\right)-\frac{1}{2}\left(\frac{a_2}{a_1}x-\frac{1}{r}x\right)\cdot\left\{-\frac{1}{ra_2}x-\left(-\frac{1}{a_1}x\right)\right\}\\
		&&-\frac{1}{2}\left(-\frac{1}{a_1}x+\frac{1}{a_1}-\frac{1}{a_1}\right)\cdot\left\{-\left(\frac{a_2}{a_1}x-\frac{a_2}{a_1}+1\right)\right\}\\
		&=&-\frac{1}{ra_1}x-\frac{1}{2a_1^2a_2r^2}(a_1-ra_2)x^2-\frac{1}{2a_1^2a_2}(a_2x+a_1-a_2)^2.
	\end{eqnarray*}
We denote this area by $V_8$.

If $-(a_1-a_2)/a_2<x<0$,
then $S(x)$ is a pentagon as in the figure \ref{fig415},
which is the same shape as that in the case where 
$-ra_2/(a_1-ra_2)<x<0$ on (1-1-1).
Hence $\area(S)=V_3$.

Therefore
	\begin{eqnarray*}
		s(X,D)&=&\vol(P_{\Sigma,D})\\
		&=&\left(\int_{-2}^{-\frac{a_1}{a_2}}V_6\,dx
		+\int_{-\frac{a_1}{a_2}}^{-\frac{a_2}{a_1-a_2}}V_5\,dx\right.\\
		&&+\left.\left.\int_{-\frac{a_2}{a_1-a_2}}^{-1}V_7\,dx
		+\int_{-1}^{-\frac{a_1-a_2}{a_2}}V_8\,dx
		+\int_{-\frac{a_1-a_2}{a_2}}^{0}V_3\,dx\right)\right|_{r=1}\\
		&=&\frac{1}{6a_1^2a_2^4(a_2-a_1)}(2a_2^6+4a_1a_2^5+2a_2^5+3a_1^2a_2^4-7a_1a_2^4-18a_1^3a_2^3+11a_1^2a_2^3\\
		&&+15a_1^4a_2^2-10a_1^3a_2^2-6a_1^5a_2+5a_1^4a_2+a_1^6-a_1^5),
	\end{eqnarray*}
(see Section \ref{maxima} (\%o7)).

\noindent(1-4-2) Suppose that $2a_1\leq (1+\sqrt{5})a_2$.
Then we have
	\[
	-(r+1)=-2<-\frac{a_2}{a_1-a_2}<-\frac{a_1}{a_2}<-1<-\frac{a_1-a_2}{a_2}<0.
	\]

If $-(r+1)=-2<x<-a_2/(a_1-a_2)$,
then $S(x)$ is a rectangle as in the figure \ref{fig421},
which is the same shape as that in the case where $-(r+1)<x<-a_1/a_2$ on (1-3).
Hence $\area(S(x))=V_6$.

If $-a_2/(a_1-a_2)<x<-a_1/a_2$,
then $S(x)$ is a pentagon as in the figure \ref{fig422}.
Hence
	\begin{eqnarray*}
		\area(S(x))&=&\left\{\frac{1}{r}x+\frac{1}{r}-(-1)\right\}\\
		&&-\frac{1}{2}\left\{\frac{1}{r}x+\frac{1}{r}-\left(\frac{a_2}{a_1}x-\frac{a_2}{a_1}+1\right)\right\}\cdot\left\{-\frac{1}{a_1}x+\frac{1}{a_1}-\left(\frac{1}{a_2}-\frac{1}{ra_2}x-\frac{1}{ra_2}\right)\right\}\\
		&=&\frac{1}{a_1r}(x+r+1)-\frac{1}{2a_1^2a_2r^2}\left\{(a_1-ra_2)x+a_1+ra_2-ra_1\right\}^2.
	\end{eqnarray*}
We denote this area by $V_9$.

If $-a_1/a_2<x<-1$,
then $S(x)$ is a hexagon as in the figure \ref{fig423},
which is the same shape as that in the case where 
$-a_2/(a_1-a_2)<x<-1$ on (1-4-1).
Hence $\area(S(x))=V_7$.

If $-1<x<0$,
then $S(x)$ is the same as that in the case where 
$-1<x<0$ on (1-4-1).

Therefore
	\begin{eqnarray*}
		s(X,D)&=&\vol(P_{\Sigma,D})\\
		&=&\left(\int_{-2}^{-\frac{a_2}{a_1-a_2}}V_6\,dx
		+\int_{-\frac{a_2}{a_1-a_2}}^{-\frac{a_1}{a_2}}V_9\,dx\right.\\
		&&+\left.\left.\int_{-\frac{a_1}{a_2}}^{-1}V_7\,dx
		+\int_{-1}^{-\frac{a_1-a_2}{a_2}}V_8\,dx
		+\int_{-\frac{a_1-a_2}{a_2}}^{0}V_3\,dx\right)\right|_{r=1}\\
		&=&\frac{1}{6a_1^2a_2^4(a_2-a_1)}(2a_2^6+4a_1a_2^5+2a_2^5+3a_1^2a_2^4-7a_1a_2^4-18a_1^3a_2^3+11a_1^2a_2^3\\
		&&+15a_1^4a_2^2-10a_1^3a_2^2-6a_1^5a_2+5a_1^4a_2+a_1^6-a_1^5),
	\end{eqnarray*}
(see Section \ref{maxima} (\%o9)).

\noindent(1-5) 
We assume that $E_t\leq K_t$.
By the same argument in the case (1-3),
we have $a_1\leq (r+1)a_2$.
We also assume that $r=1$.
Then $G_t\leq I_t$.
We also assume that $3a_2\geq 2a_1$.
Then $F_t\leq L_t$.
(See the figure \ref{fig501}.)

If $B_t=F_t$, then $-x/a_1+1/a_1=2/a_2$.
Hence $x=-(2a_1-a_2)/a_2$.
By the same argument in the case (1-3), we have the following two values of $x$.
If $A_t=E_t$, then $x=-a_1/a_2$.
If $B=K$, then $x=-r$, in particular $x=-1$.
By the same argument in the case (1-2), we have the following two values of $x$.
If $C=L$ and $D=K$,
then $x=-(r+1)$, in particular, $x=-2$.
If $B_t=E_t$, then $x=-(a_1-a_2)/a_2$.
By the same argument in the case (1-1), we have the following two values of $x$.
If $C=I$ and $D=H$,
then $x=-1$.
If $A=O$ and $B=H$,
then $x=0$.

Since $2a_1\leq 3a_2$, we have $2a_1-a_2\leq 2a_2$.
Hence $-2\leq -(2a_1-a_2)/a_2$.
Since $a_1<a_1+(a_1-a_2)$, we have $-(2a_1-a_2)/a_2<-a_1/a_2$.
Since $a_2<a_1$, we have $-a_1/a_2<-1$.
Since $a_1\leq 2a_2$, we have $a_1-a_2\leq a_2$.
Hence $-1\leq -(a_1-a_2)/a_2$.
Then we have
	\[
	-(r+1)=-2<-\frac{2a_1-a_2}{a_2}<-\frac{a_1}{a_2}<-1<-\frac{a_1-a_2}{a_2}<0.
	\]

For various values of $x$,
we consider shapes of $S(x)$ and areas of them.

If $-(r+1)=-2<x<-(2a_1-a_2)/a_2$,
then $S(x)$ is a trapezoid as in the figure \ref{fig502}.
Hence
	\begin{eqnarray*}
		\area(S(x))&=&\frac{1}{2}\left\{\left(\frac{1}{a_2}x-\frac{1}{ra_2}x-\frac{1}{ra_2}
-\left(-\frac{1}{a_1}x\right)\right)\right\}\cdot\left\{\frac{1}{r}x+\frac{1}{r}-(-1)\right\}\\
 &=&\frac{1}{2a_1r^2a_2}(x+r+1)\left\{(2a_2r-a_1)x+3a_1r-a_1\right\}.
	\end{eqnarray*}

If $-(2a_1-a_2)/a_2<x<-a_1/a_2$,
then $S(x)$ is a pentagon as in the figure \ref{fig503},
which is the same shape as that in case where 
$-a_2/(a_1-a_2)<x<-a_1/a_2$ on (1-4-2).
Hence $\area(S(x))=V_9$.

If $-a_1/a_2<x<-1$,
then $S(x)$ is a hexagon as in the figure \ref{fig504},
which is the same shape as that in case where 
$-a_2/(a_1-a_2)<x<-1$ on (1-4-1).
Hence $\area(S(x))=V_7$.

If $-1<x<-(a_1-a_2)/a_2$,
then $S(x)$ is a hexagon as in the figure \ref{fig505},
which is the same shape as that in case where 
$-1<x<-(a_1-a_2)/a_2$ on (1-4-1).
Hence $\area(S(x))=V_8$.

If $-(a_1-a_2)/a_2<x<0$,
then $S(x)$ is a pentagon as in the figure \ref{fig506},
which is the same shape as that in case where 
$-ra_2/(a_1-ra_2)<x<0$ on (1-1-1).
Hence $\area(S(x))=V_3$.

Therefore
	\begin{eqnarray*}
		s(X,D)&=&\vol(P_{\Sigma,D})\\
		&=&\left(\int_{-2}^{-\frac{2a_1-a_2}{a_2}}\frac{1}{2a_1r^2a_2}(x+r+1)\left\{(2a_2r-a_1)x+3a_1r-a_1\right\}dx\right.\\
		&&+\int_{-\frac{2a_1-a_2}{a_2}}^{-\frac{a_1}{a_2}}V_9\,dx
		+\left.\left.\int_{-\frac{a_1}{a_2}}^{-1}V_7\,dx
		+\int_{-1}^{-\frac{a_1-a_2}{a_2}}V_8\,dx
		+\int_{-\frac{a_1-a_2}{a_2}}^{0}V_3\,dx\right)\right|_{r=1}\\
		&=&\frac{1}{6a_1^2a_2^4}(2a_2^5-21a_1a_2^4+2a_2^4+36a_1^2a_2^3-5a_1a_2^3-18a_1^3a_2^2\\
		&&+6a_1^2a_2^2+5a_1^4a_2-4a_1^3a_2-a_2^5+a_1^4),
	\end{eqnarray*}
(see Section \ref{maxima} (\%o11)).


\noindent(2) Suppose that $P_D$ is a triangle, i.e., $a_1\leq ra_2$. Then the facets 
of $P_D$ are defined by the equations $x=-a_1$, $y=0$, and $-x+ry=0$. Hence we 
have $I_D=\{1,3,4\}$ and 
$$P_{\Sigma, D}=
\left\{(x, y, t)\in (M_D\times \ZZ)_{\R}\left|
\begin{array}{l}
0\leq x+a_1 t <1, \\
0\leq -y<1, \\
0\leq -x+ry<1
\end{array}
\right.\right\}.$$

The point of the intersection of the three hyperplanes defined by $x+a_1t=0$, $-y=0$, 
and $-x+ry=1$ (resp. $x+a_1t=0$, $-y=1$, and $-x+ry=0$; $x+a_1t=1$, $-y=0$, and 
$-x+ry=0$) is $(-1, 0, 1/a_1)$ (resp. $(-r, -1, r/a_1)$; $(0, 0, 1/a_1)$). 

Let $v_1={}^t(-1, 0, 1/a_1)$, $v_2={}^t(-r, -1, r/a_1)$, and $v_3={}^t(0, 0, 1/a_1)$. Since 
$P_{\Sigma, D}$ is spanned by $v_1$, $v_2$, and $v_3$, we have 
$$s(X,D)=\vol (P_{\Sigma, D})= |\mathrm{det}(v_1,v_2,v_3)|=\dfrac{1}{a_1}.$$

\noindent(3) Suppose that $P_D$ is not full-dimensional. $P_D$ is not full-dimensional 
if and only if $a_1=0$ or $a_2=0$. First we consider the case where $a_2=0$. Then 
$P_D$ is a line segment whose edges are defined by $x=-a_1$ and $-x+ry=0$. Hence 
we have $I_D=\{1,4\}$. We see that $M_D =\ZZ\times \{0\}\subset M$ with basis 
$\{\langle1, 0\rangle\}$. Let $\pi :N\rightarrow N_D=\ZZ\times \{0\}$ be the first 
projection map. Since $v_1\cdot \langle1, 0\rangle=1$ (resp. 
$v_4\cdot \langle1, 0\rangle=-1$), $\pi (v_1)$ (resp. $\pi (v_4)$) is $N_D$-primitive 
and $c_1=1$ (resp. $c_4=1$). Therefore we have 
$${P'}_{\Sigma, D}=
\left\{(x , t)\in (M_D\times \ZZ)_{\R}\left|
\begin{array}{l}
0\leq x+a_1 t <1, \\
0\leq -x<1
\end{array}
\right.\right\}.$$
${P'}_{\Sigma, D}$ is the parallelogram spanned by the two vectors $v_1={}^t(0, 1/a_1)$ 
and $v_2={}^t(-1,1/a_1)$. Hence we have 
$$s(X, D)=\vol ({P'}_{\Sigma, D})= |\mathrm{det}(v_1,v_2)|=\dfrac{1}{a_1}.$$

Next we consider the case where $a_1=0$. Since $\mathrm{Sec} (X,a_2D_2)\cong k[x]$, 
the $F$-signature is not defined. (An $F$-signature is defined for the section ring 
whose dimension is at least two.)
\end{exa}

\section{Figures}\label{fig}
In this section, we gather together the figures in the computation of $F$-signatures of 
Hirzebruch surfaces. 

\clearpage
\begin{figure}[H]
\centering\scalebox{0.8}{
{\unitlength 0.1in%
}%
}\caption{}\label{fig001}
\centering\scalebox{0.8}{
{\unitlength 0.1in%
%
}%
}\caption{}\label{fig002}
\end{figure}
\begin{figure}[H]
\centering\scalebox{0.8}{
{\unitlength 0.1in%
%
}%
}\caption{}\label{fig003}
\centering\scalebox{0.8}{
{\unitlength 0.1in%
%
}%
}\caption{}\label{fig101}
\end{figure}
\begin{figure}[H]
\centering\scalebox{0.8}{
{\unitlength 0.1in%
%
}%
}\caption{}\label{fig102}
\centering\scalebox{0.8}{
{\unitlength 0.1in%
%
}%
}\caption{}\label{fig111}
\end{figure}
\begin{figure}[H]
\centering\scalebox{0.8}{
{\unitlength 0.1in%
%
}%
}\caption{}\label{fig506}
\end{figure}
\section{Maxima source codes}\label{maxima}
To compute integrals, we use the Maxima which is a computer algebra system.
In particular, we use a software for Maxima called wxMaxima.
In the following code,
we use $a$ (resp. $b$) instead of $a_1$ (resp. $a_2$)
The following are source code:

\begin{verbatim}
(%i1)/*case 1 (2a_2\leq a_1)*/
factor(integrate((1/(2*a^2*b*r^2))*((a-b*r)*x+a+b*r)^2,
x,(b*r+a)/(b*r-a),a/(b*r-a))
+integrate((1/a)*(1/r-b/a)*x+(1/(2*a))*(2/r+b/a),x,a/(b*r-a),-1)
+integrate((b/(2*a^2))*(1-2*x),x,-1,(b*r)/(b*r-a))
+integrate(-(x/(a*r))-(1/2)*(b/a-1/r)*(1/a-1/(b*r))*x^2,
x,(b*r)/(b*r-a),0));
\end{verbatim}

$$-{{b}\over{a\,\left(b\,r-a\right)}}\leqno{\tt (\%o1)}$$

\begin{verbatim}
(%i2)/*case 1 (2a_2> a_1)*/
factor(integrate(1/(2*a^2*b*r^2)*((a-r*b)*x+a+r*b)^2,
x,-(a+r*b)/(a-b*r),-a/(a-r*b))
+integrate(1/(2*r*a^2)*(2*(a-r*b)*x+2*a+r*b),
x,-a/(a-r*b),-(r*b)/(a-r*b))
+integrate(1/(r*a)-(r*b-a)^2/(2*r^2*a^2*b)*x^2,
x,-(r*b)/(a-r*b),-1)
+integrate(-x/(r*a)-(r*b-a)^2*x^2/(2*r^2*a^2*b),x,-1,0));
\end{verbatim}

$$-{{b}\over{a\,\left(b\,r-a\right)}}\leqno{\tt (\%o2)}$$

\begin{verbatim}
(%i3)/*case_2 r=1 */
factor(integrate(1/(2*a*b*r^2)*((a-2*b*r)*x+2*b*r+a-a*r)*(x+r+1),
x,-(r+1),-a/(a-r*b))
+integrate((x+r+1)/(r*a)-(b*x+a)^2/(2*a^2*b),x,-a/(a-r*b),-(a-b)/b)
+integrate(1/(2*a^2*r)*(2*(a-r*b)*x+2*a+r*b),x,-(a-b)/b,-1)
+integrate(b*(1-2*x)/(2*a^2),x,-1,-(r*b)/(a-r*b))
+integrate(-x/(r*a)-(r*b-a)^2*x^2/(2*r^2*a^2*b),x,-r*b/(a-r*b),0));
\end{verbatim}

$${{b^3\,r^3+6\,b^3\,r^2-3\,a\,b^2\,r^2+12\,b^3\,r-12\,a\,b^2\,r+3\,a
 ^2\,b\,r+2\,b^3-12\,a\,b^2+6\,a^2\,b-a^3}\over{6\,a\,b^2\,\left(b\,r
 -a\right)}}\leqno{\tt (\%o3)}$$
 
 \begin{verbatim}
(%i4)/*case_2 r>1 */
factor(integrate(1/(2*a*b*r^2)*((a-2*b*r)*x+2*b*r+a-a*r)*(x+r+1),
x,-(r+1),-a/(a-r*b))
+integrate((x+r+1)/(r*a)-(b*x+a)^2/(2*a^2*b),x,-a/(a-r*b),-(a-b)/b)
+integrate(1/(2*a^2*r)*(2*(a-r*b)*x+2*a+r*b),x,-(a-b)/b,-r*b/(a-r*b))
+integrate(1/(r*a)-(r*b-a)^2*x^2/(2*r^2*a^2*b),x,-(r*b)/(a-r*b),-1)
+integrate(-x/(r*a)-(r*b-a)^2*x^2/(2*r^2*a^2*b),x,-1,0));
\end{verbatim}

$${{b^3\,r^3+6\,b^3\,r^2-3\,a\,b^2\,r^2+12\,b^3\,r-12\,a\,b^2\,r+3\,a
 ^2\,b\,r+2\,b^3-12\,a\,b^2+6\,a^2\,b-a^3}\over{6\,a\,b^2\,\left(b\,r
 -a\right)}}\leqno{\tt (\%o4)}$$
 
 \begin{verbatim}
(%i5)/*case3 */
factor(integrate((x+r+1)/(a*r),x,-(r+1),-a/b)
+integrate((x+r+1)/(r*a)-(b*x+a)^2/(2*a^2*b),x,-a/b,-r)
+integrate(1/(r*a)-(r*b-a)^2*x^2/(2*r^2*a^2*b),x,-r,-1)
+integrate(-x/(r*a)-(r*b-a)^2*x^2/(2*r^2*a^2*b),x,-1,0));
\end{verbatim}

$$-{{b^2\,r^2-2\,a\,b\,r-6\,b^2+a^2}\over{6\,a\,b^2}}\leqno{\tt (\%o5)}$$

\begin{verbatim}
(%i6)/* case4 2a_1>(1+\sqrt{5})*a_2 */
f(r):=integrate((x+r+1)/(a*r),x,-2,-a/b)
+integrate((x+r+1)/(r*a)-(b*x+a)^2/(2*a^2*b),x,-a/b,-b/(a-b))
+integrate((x+r+1)/(r*a)-(b*x+a)^2/(2*a^2*b)
-((a-r*b)*x+a+r*b-r*a)^2/(2*a^2*b*r^2),x,-b/(a-b),-1)
+integrate(-x/(r*a)-(a-r*b)*x^2/(2*a^2*b*r^2)-(b*x+a-b)^2/(2*a^2*b),
x,-1,-(a-b)/b)
+integrate(-x/(r*a)-(r*b-a)^2*x^2/(2*r^2*a^2*b),x,-(a-b)/b,0);
\end{verbatim}

\begin{verbatim}
(%i7)factor(f(1));
\end{verbatim}
{\tiny{
$${{2\,b^6+4\,a\,b^5+2\,b^5+3\,a^2\,b^4-7\,a\,b^4-18\,a^3\,b^3+11\,a^
 2\,b^3+15\,a^4\,b^2-10\,a^3\,b^2-6\,a^5\,b+5\,a^4\,b+a^6-a^5}\over{6
 \,a^2\,b^4\,\left(b-a\right)}}\leqno{\tt (\%o7)}$$}}
 
 \begin{verbatim}
(%i8)/* case4 2a_1\leq(1+\sqrt{5})*a_2 */
g(r):=integrate((x+r+1)/(a*r),x,-2,-b/(a-b))
+integrate((x+r+1)/(r*a)-((a-r*b)*x+a+r*b-r*a)^2/(2*a^2*b*r^2),
x,-b/(a-b),-a/b)
+integrate((x+r+1)/(r*a)-(b*x+a)^2/(2*a^2*b)
-((a-r*b)*x+a+r*b-r*a)^2/(2*a^2*b*r^2),x,-a/b,-1)
+integrate(-x/(r*a)-(a-r*b)*x^2/(2*a^2*b*r^2)-(b*x+a-b)^2/(2*a^2*b),
x,-1,-(a-b)/b)
+integrate(-x/(r*a)-(r*b-a)^2*x^2/(2*r^2*a^2*b),x,-(a-b)/b,0);
\end{verbatim}

\begin{verbatim}
(%i9)factor(g(1));
\end{verbatim}

{\tiny{$${{2\,b^6+4\,a\,b^5+2\,b^5+3\,a^2\,b^4-7\,a\,b^4-18\,a^3\,b^3+11\,a^
 2\,b^3+15\,a^4\,b^2-10\,a^3\,b^2-6\,a^5\,b+5\,a^4\,b+a^6-a^5}\over{6
 \,a^2\,b^4\,\left(b-a\right)}}\leqno{\tt (\%o9)}$$}}
 
 \begin{verbatim}
(%i10) /*case 5*/
h(r):=integrate((x+r+1)*((2*b*r-a)*x+3*a*r-a)/(2*a*r^2*b),
x,-2,-(2*a-b)/b)
+integrate((x+r+1)/(r*a)-((a-r*b)*x+a+r*b-r*a)^2/(2*a^2*b*r^2),
x,-(2*a-b)/b,-a/b)
+integrate((x+r+1)/(r*a)-(b*x+a)^2/(2*a^2*b)
-((a-r*b)*x+a+r*b-r*a)^2/(2*a^2*b*r^2),x,-a/b,-1)
+integrate(-x/(r*a)-(a-r*b)*x^2/(2*a^2*b*r^2)-(b*x+a-b)^2/(2*a^2*b),
x,-1,-(a-b)/b)
+integrate(-x/(r*a)-((r*b-a)*x)^2/(2*r^2*a^2*b),x,-(a-b)/b,0);
\end{verbatim}

\begin{verbatim}
(%i11)factor(h(1));
\end{verbatim}

$${{2\,b^5-21\,a\,b^4+2\,b^4+36\,a^2\,b^3-5\,a\,b^3-18\,a^3\,b^2+6\,a
 ^2\,b^2+5\,a^4\,b-4\,a^3\,b-a^5+a^4}\over{6\,a^2\,b^4}}\leqno{\tt (\%o11)}$$


\begin{thebibliography}{1}

\bibitem[CLS11]{CLS}
David~A. Cox, John~B. Little, and Henry~K. Schenck, \emph{Toric varieties},
  Graduate Studies in Mathematics, vol. 124, American Mathematical Society,
  Providence, RI, 2011. \MR{2810322}

\bibitem[HN15]{HN}
Mitsuyasu Hashimoto and Yusuke Nakajima, \emph{Generalized {$F$}-signature of
  invariant subrings}, J. Algebra \textbf{443} (2015), 142--152. \MR{3400399}

\bibitem[HL02]{HL}
Craig Huneke and Graham~J. Leuschke, \emph{Two theorems about maximal
  {C}ohen-{M}acaulay modules}, Math. Ann. \textbf{324} (2002), no.~2, 391--404.
  \MR{1933863}

\bibitem[San15]{San}
Akiyoshi Sannai, \emph{On dual {$F$}-signature}, Int. Math. Res. Not. IMRN
  (2015), no.~1, 197--211. \MR{3340299}

\bibitem[Tuc12]{T}
Kevin Tucker, \emph{{$F$}-signature exists}, Invent. Math. \textbf{190} (2012),
  no.~3, 743--765. \MR{2995185}

\bibitem[Kor12]{K}
Michael~R. Von~Korff, \emph{The {F}-{S}ignature of {T}oric {V}arieties},
  ProQuest LLC, Ann Arbor, MI, 2012, Thesis (Ph.D.)--University of Michigan.
  \MR{3093997}

\end{thebibliography}

\end{document}